\documentclass[a4wide,11pt]{article}
\usepackage{amsmath,amsfonts,amssymb,rotate,lscape}
\usepackage{epsfig}
\usepackage{a4wide,graphicx}
\usepackage{epstopdf}
\usepackage[pdftex]{rotating}
\usepackage[latin1]{inputenc}
\usepackage[francais,english]{babel}
\usepackage[T1]{fontenc}
\usepackage{natbib}
\usepackage{color}
\definecolor{linkcol}{rgb}{0,0,0.4}
\definecolor{citecol}{rgb}{0.5,0,0}
\usepackage{xcolor}
\usepackage{colortbl,hhline}
\definecolor{lightgray}{gray}{0.90}

\batchmode

\begin{document}

\title{A Generalization of the Exponential-Logarithmic Distribution for Reliability and Life Data Analysis}

\author{Mohieddine Rahmouni\thanks{\textbf{Corresponding author}, ESSECT, University of Tunis, Tunisia, and Community College in Abqaiq, King Faisal University, Saudi Arabia. E-mail address: \texttt{mohieddine.rahmouni@gmail.com}}, \ Ayman Orabi\thanks{Community College in Abqaiq, King Faisal University, Saudi Arabia.}\\
\\}

\date{January 5,2018}
\maketitle

\begin{abstract}

In this paper, we introduce a new two-parameter lifetime distribution, called the exponential-generalized truncated logarithmic (EGTL) distribution, by compounding the exponential and generalized truncated logarithmic distributions. Our procedure generalizes the exponential-logarithmic (EL) distribution modelling the reliability of systems by the use of first-order concepts, where the minimum lifetime is considered \citep{Tahmasbi2008}. In our approach, we assume that a system fails if a given number $k$ of the components fails and then, we consider the $k^{th}$-smallest value of lifetime instead of the minimum lifetime.
The reliability and failure rate functions as well as their properties are presented for some special cases. The estimation of the parameters is attained by the maximum likelihood, the expectation maximization algorithm, the method of moments and the Bayesian approach, with a simulation study performed to illustrate the different methods of estimation. The application study is illustrated based on two real data sets used in many applications of reliability.\\


\end{abstract}

\section{Introduction}

Lifetime distributions are often used in reliability theory and survival analysis for modelling real data. They play a fundamental role in reliability in diverse disciplines such as finance, manufacture, biological sciences, physics and engineering. The exponential distribution is a basic model in reliability theory and
survival analysis. It is often used to model system reliability at a component level, assuming the failure rate is constant \citep{Balakrishnan1995,Barlow1975,Sinha1980}. In recent years, a growing number of scholarly papers has been devoted to accommodate lifetime distributions with increasing or decreasing failure rate functions. The motivation is to give a parametric fit for real data sets where the underlying failure rates, arising on a latent competing risk problem base, present monotone shapes (nonconstant hazard rates).
The proposed distributions are introduced as extensions of the exponential distribution, following \cite{Adamidis1998} and \cite{Kus2007}, by compounding some useful lifetime and truncated discrete distributions (for review, see \cite{Barreto-Souza2009,Chahkandi2009,Silva2010,Barreto-Souza2011,Cancho2011,Louzada-Neto2011,Morais2011,Hemmati2011,Nadarajah2013,Bakouch2014}, and others).
The genesis is stated on competing risk scenarios in presence of latent
risks, i.e. there is no information about the causes of the component's failure \citep{Basu1982}. In fact, a system  may  experience  multiple  failure  processes  that  compete  against  each other, and whichever occurs first  can  cause the  system  to fail \citep{Rafiee2017,Kalbfleisch2002,Andersen2002,Tsiatis1998}. The term competing risks refers to duration data where two or more causes are competing to determine the observed time-to-failure. The potential multiple causes of failure are not mutually exclusive but the interest lies in the time to the first coming one \citep{Putter2007,Bakoyannis2012}. For further details see \citet{Basu1981}.

In the same way, the exponential-logarithmic (EL) distribution was proposed by \cite{Tahmasbi2008} as a log-series mixture of exponential random variables. This two-parameter distribution with decreasing failure rate (DFR) is obtained by mixing the exponential and logarithmic distributions.
It is based on the idea of modelling the system's reliability where the time-to-failure occurs due to the presence of an unknown number of initial defects of some components is considered. Suppose the breakdown (failure) of a system of components occurs due to the presence of a non-observable number, $Z$, of initial defects of the same kind, that can be identifiable only after causing failure and are repaired perfectly \citep{Adamidis1998,Kus2007}. Let $T_{i}$ be the failure time of the system due to the $i^{th}$ defect, for $i \geq 1$. If we assume that $T=(T_{1}, T_{2}, ..., T_{Z})$ are iid exponential random variables independent of $Z$, that follows a truncated logarithmic distribution,
then the time to the first failure is adequately modelled by the
EL distribution \citep{Barreto-Souza2015,Bourguignon2014,Ross1976}.
For reliability studies, $X_{(1)}=min\{T_i\}_{i=1}^{Z}$ and $X_{(Z)}=max\{T_i\}_{i=1}^{Z}$ are used respectively in serial and parallel systems with identical components \citep{Chahkandi2009,Ramos2015}.
However, one may determine the distribution of the $k^{th}$ smallest value of the time-to-failure ($k^{th}$ order
statistic), instead of the minimum lifetime (first order statistic).

There is a huge literature on the order
statistics for reliability engineering (for review, see \cite{Barlow1975,Sarhan1962,Barlow1965,pyke1965,Gnedenko1969,Pledger1971,Barlow1981,David1981,Bain1991}, and references contained therein).
The motivation arises in
reliability theory, where the so-called $k$-out-of-$n$ systems
are studied \citep{Xie2008,Xie2005,Proschan1976,Kim1988}.
An engineering system consisting of $n$ components is working if at least $k$ out of the total $n$ components are operating and it breaks down if $(n-k+1)$ or more components fail. Hence, a $k$-out-of-$n$ system fails at the time of the $(n-k+1)^{th}$ component failure \citep{Barlow1975,Kamps1995,Cramer2001}. This binary-state context is based on the assumption that a system or its components can be either fully working or completely failing.
However, in reality, a system may provide its specified function at less than full capacity when
some of its components operate in a degraded state \citep{Ramirez-Marquez2005}.
The binary $k$-out-of-$n$ system reliability models have been extended to multi-state $k$-out-of-$n$ system reliability models by allowing more than two performance levels for the system and its components \citep{Eryilmaz2014}. Multi-state systems contain units presenting two or more failed states with multiple modes of failure and one working state \citep{Anzanello2009}.
Reliability models provide, through multi-state, more realistic representations of engineering systems \citep{Yingkui2012}.
Many authors have made contributions about the reliability estimation approaches for multi-state systems \citep{Jenney1986,Page1988,Rocco2005,Ramirez-Marquez2004,Ramirez-Marquez2008,Levitin2007,Levitin2008}.

In this paper, we generalize the EL distribution \citep{Tahmasbi2008} modelling the time to the first failure, to a
distribution more appropriate for modelling any order statistic (second, third, or any  $k^{th}$ lifetime). For instance, suppose a machine produces a random number, $Z$ units, of light bulbs or wire fuses which are put through a life test. Each item has a random lifetime $T_{i}$, $i=1, 2, ...,Z$.
The EL distribution \citep{Tahmasbi2008} is focused only on the minimum time-to-failure of the first of the $Z$ functioning components.
However, we may be interested in the $k^{th}$ duration and then determine the lifetime distribution for the order statistics, assuming the system will fail if $k$ of the units fail. We may let $X_{(1)}  < X_{(2)}  < ... <X_{(Z)}$  be the order statistics of $z$ independent observations of the time $T_{i}$ and then, we consider the $k^{th}$-smallest value of lifetime instead of the minimum lifetime.
We assume that the $T_{i}$'s are not observable, but that $X_{(k)}$ is. We would like to estimate the lifetime distribution given the observations on $X_{(k)}$.
The proposed new family of lifetime distributions is obtained by compounding the exponential and generalized truncated logarithmic distributions, named exponential-generalized truncated logarithmic (EGTL) distribution.

The paper is organized as follows: In section 2, we present the new family of lifetime distributions and the probability density function (pdf) for some special cases. The moment generating function, the $r^{th}$ moment, the reliability, the failure rate function and the random number generation are discussed in this section. The estimation of parameters for this new family of   distributions will be discussed in section 3.
It is attained by maximum likelihood (MLEs) and expectation maximization (EM) algorithms. The method of moments and the Bayesian approach are also presented as possible alternatives to the MLEs method. As illustration of these three methods of estimation, numerical computations will be performed in section 4.
The application study is illustrated based on two real data sets in section 5. The last section concludes the paper.

\section{Properties of the distribution}
\subsection{Distribution}
The derivation of the new family of lifetime distributions depends on the generalization of the compound exponential and truncated logarithmic distributions as follows:
Let $T=(T_{1}, T_{2}, ..., T_{Z})$ be iid exponential random variables with scale parameter $\theta >0$ and a pdf given by: $f(t)=\theta e^{-\theta t}$ , for $t>0$, where $Z$  is a discrete random variable following a logarithmic-series distribution with parameter $0<p<1$  and a probability mass function (pmf), $P(Z=z)$, given by:
\begin{equation}\label{eq:logarithmic}
   P(Z=z) = \frac{1}{-\ln(1-p)} \frac{p^z}{z} ; z \in \{1,2,3,\dots\}
\end{equation}

If $Z$ is a truncated at $k-1$ logarithmic random variable with parameter $p$, then the probability function $P_{k-1}(Z=z)$ will be given by:
\begin{equation}\label{eq:tr-logarithmic}
   P_{k-1}(Z=z) = \frac{1}{A(p,k)} \frac{p^{z}}{z} ; k = 1,2,3,\dots\ , z \mbox{ and } z=k, k+1,\dots\
\end{equation}

where,
\begin{equation}\label{eq:tr-logarithmic-A}
   A(p,k) = \sum_{j=k}^{\infty}\frac{p^{j}}{j}
          = - \ln(1-p)-\psi(k)\sum_{j=1}^{k-1}\frac{p^{j}}{j}
\end{equation}

and,
\begin{equation}\label{eq:psi-k}
   \psi(k) =
   \begin{cases}
   0 & \mbox{if } k=1 \\
   1 & \mbox{if } k=2, 3, ...,z
   \end{cases}
\end{equation}

The pdf of the $k^{th}$  order statistic $X_{(k)}$ (the $k^{th}$-smallest value of lifetime)  exponentially distributed is given by the equation (\ref{eq:order}) (see, \cite{David1970,Balakrishnan1991,Balakrishnan1996}):

\begin{equation}\label{eq:order}
   f_{k}(x/z,\theta)=\frac{\theta \Gamma(z+1)}{\Gamma(k)\Gamma(z-k+1)}e^{-\theta(z-k+1)x}(1-e^{-\theta x})^{k-1} \mbox{ ; } \theta, x > 0
\end{equation}

From equations  (\ref{eq:tr-logarithmic}) and (\ref{eq:order}) the joint probability density is derived as\footnote{The proofs of all steps and equations are presented in the appendix.}:

\begin{equation}\label{eq:joint-dist}
   g_{k}(x,z/p,\theta)=\frac{\Gamma(z)}{\Gamma(k)\Gamma(z-k+1)}\frac{\theta p^{z} e^{-\theta (z-k+1) x} (1-e^{-\theta x})^{k-1}}{A(p,k)}
\end{equation}

where, $x$ is the lifetime of a system and $z$ is the last order statistic. In equation (\ref{eq:joint-dist}) we consider the ascending order $X_{(1)}  < X_{(2)}  < ... <X_{(Z)}$. The joint probability density is determined by compounding a truncated at $k-1$  logarithmic  series distribution and the pdf of the $k^{th}$  order statistic ($k=1, 2, ..., z$). The use of the truncated at $k-1$  logarithmic distribution is motivated by mathematical interest because we are interested in the $k^{th}$ order statistic. There is a left-truncation scheme, where only ($z-k+1$) individuals who survive a sufficient time are included, i.e. we observe only individuals or units with $X_{(k)}$ exceeding the time of the event that truncates individuals. In comparison with the formulation of \cite{Tahmasbi2008} and \cite{Adamidis1998}, we consider the $k^{th}$-smallest value of lifetime instead of the minimum lifetime $X_{(1)}=min\{T_i\}_{i=1}^{Z}$.
So, our proposed new lifetime distribution, named exponential-generalized truncated logarithmic (EGTL) distribution, is the marginal density distribution of $x$ given by:

\begin{equation}\label{eq:marginal-x-dist}
   g_{k}(x/p,\theta, k)=\frac{\theta p^{k} e^{-\theta x} (1-e^{-\theta x})^{k-1}}{A(p,k)(1-pe^{-\theta x})^{k}}
    \mbox{ ; }  \quad x \in [0, \infty)
\end{equation}

where $0<p<1$ is the shape parameter and $\theta$ is the scale parameter.
This distribution is more appropriate for modelling any $k^{th}$ order statistic ($2^{nd}$, $3^{rd}$ or any $k^{th}$ lifetime).
The particular case of the EGTL density function, for $k=1$, is the EL distribution modelling the time of the first failure, $X_{(1)}=min\{T_i\}_{i=1}^{Z}$, given by \cite{Tahmasbi2008}:

\[g_{1}(x; p, \theta) = \frac{p\theta e^{-\theta x}}{-\ln (1-p) (1-pe^{-\theta x})}\]

\begin{figure}[htp]
\begin{center}
\begin{tabular}[c]{cc}
\fbox{\includegraphics[scale=0.23]{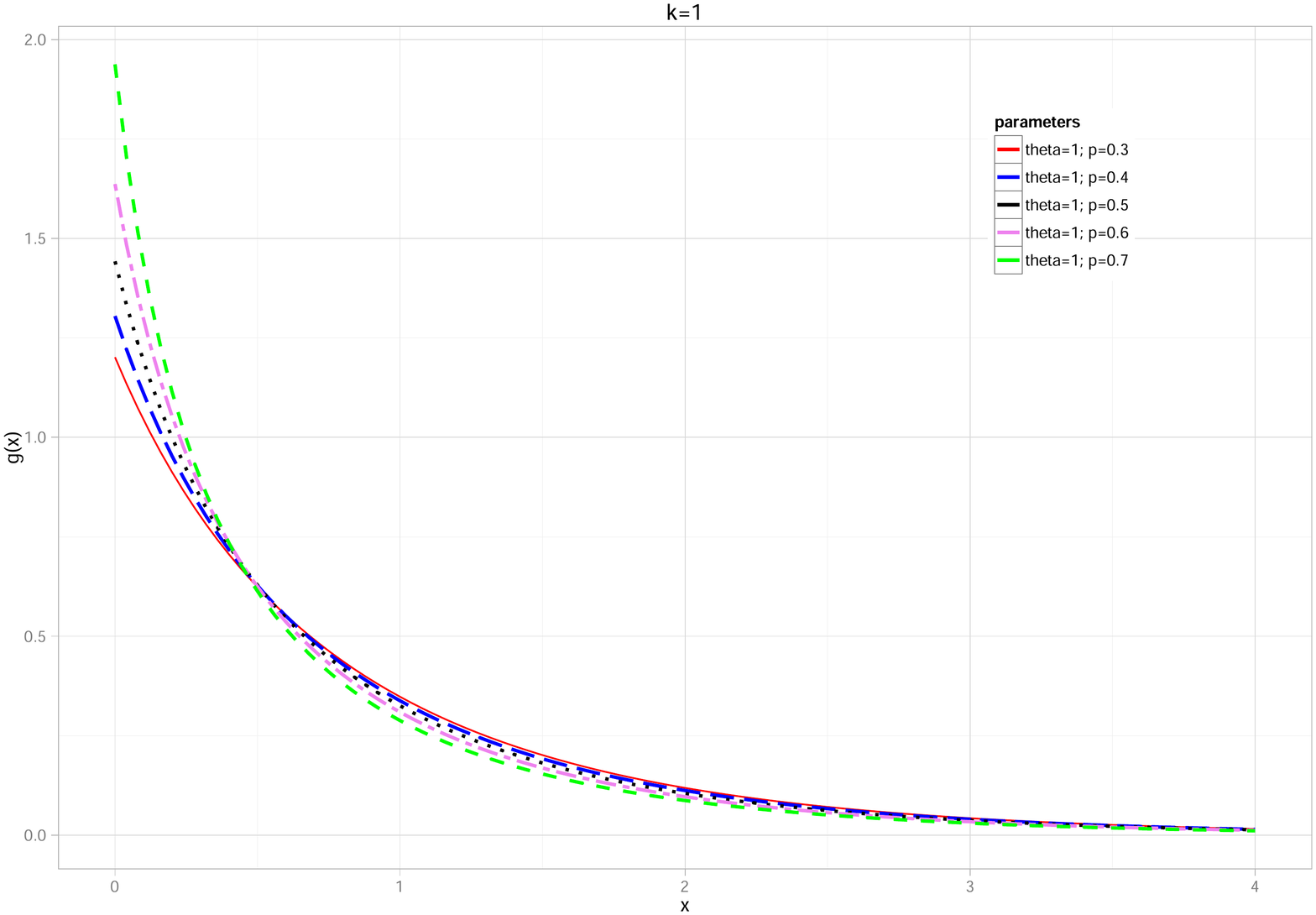}} &
\fbox{\includegraphics[scale=0.23]{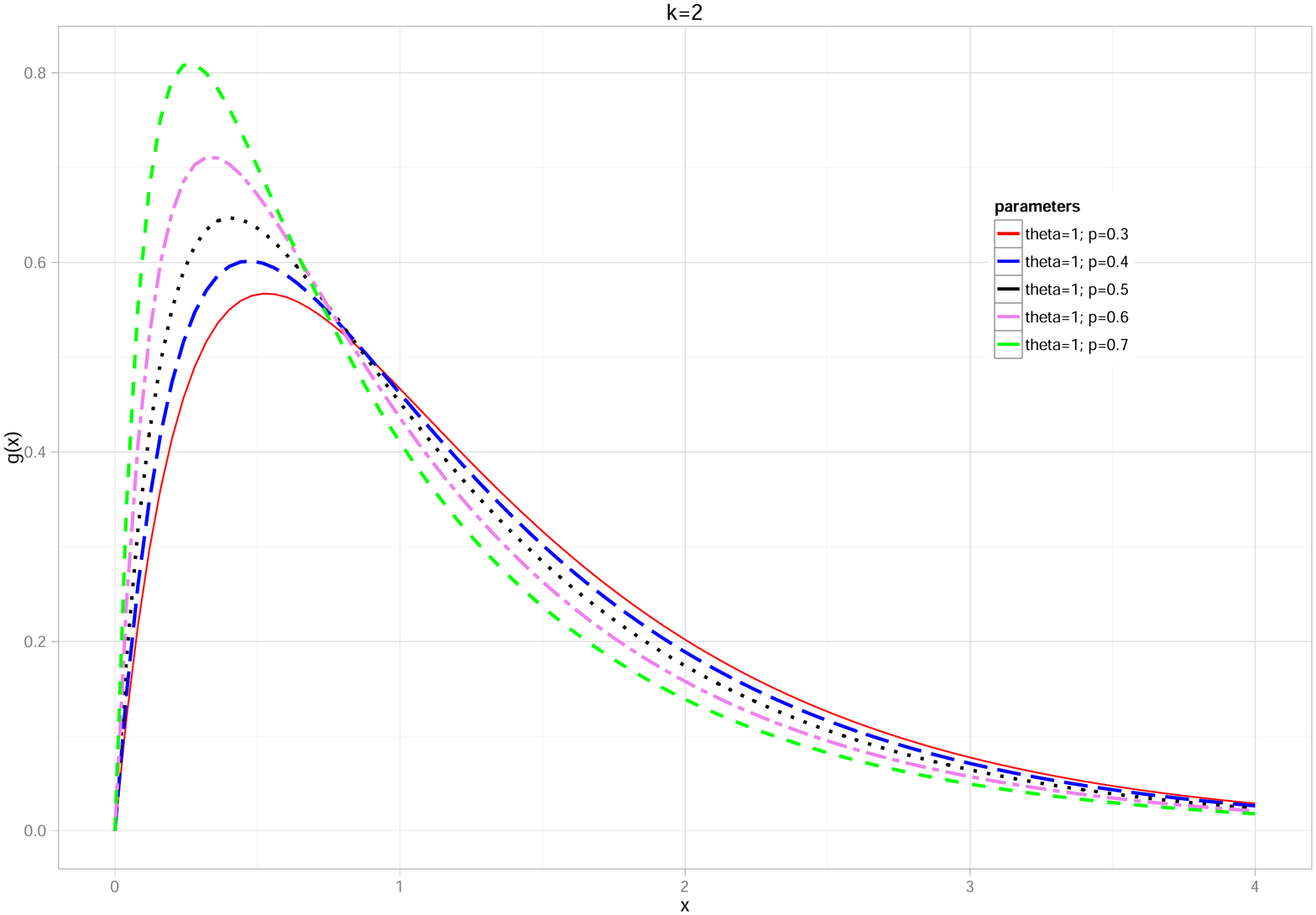}} \\
\fbox{\includegraphics[scale=0.23]{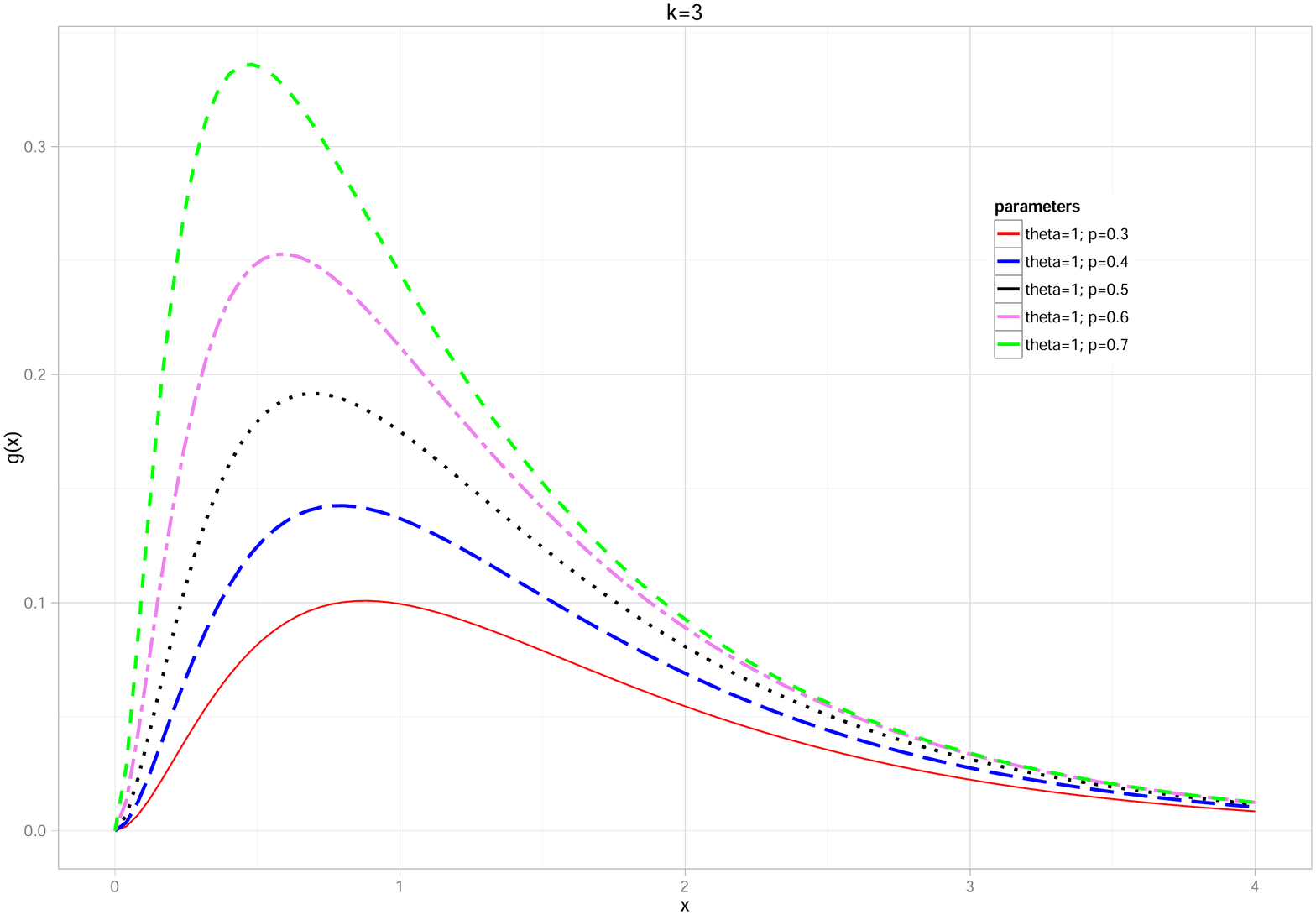}} &
\fbox{\includegraphics[scale=0.23]{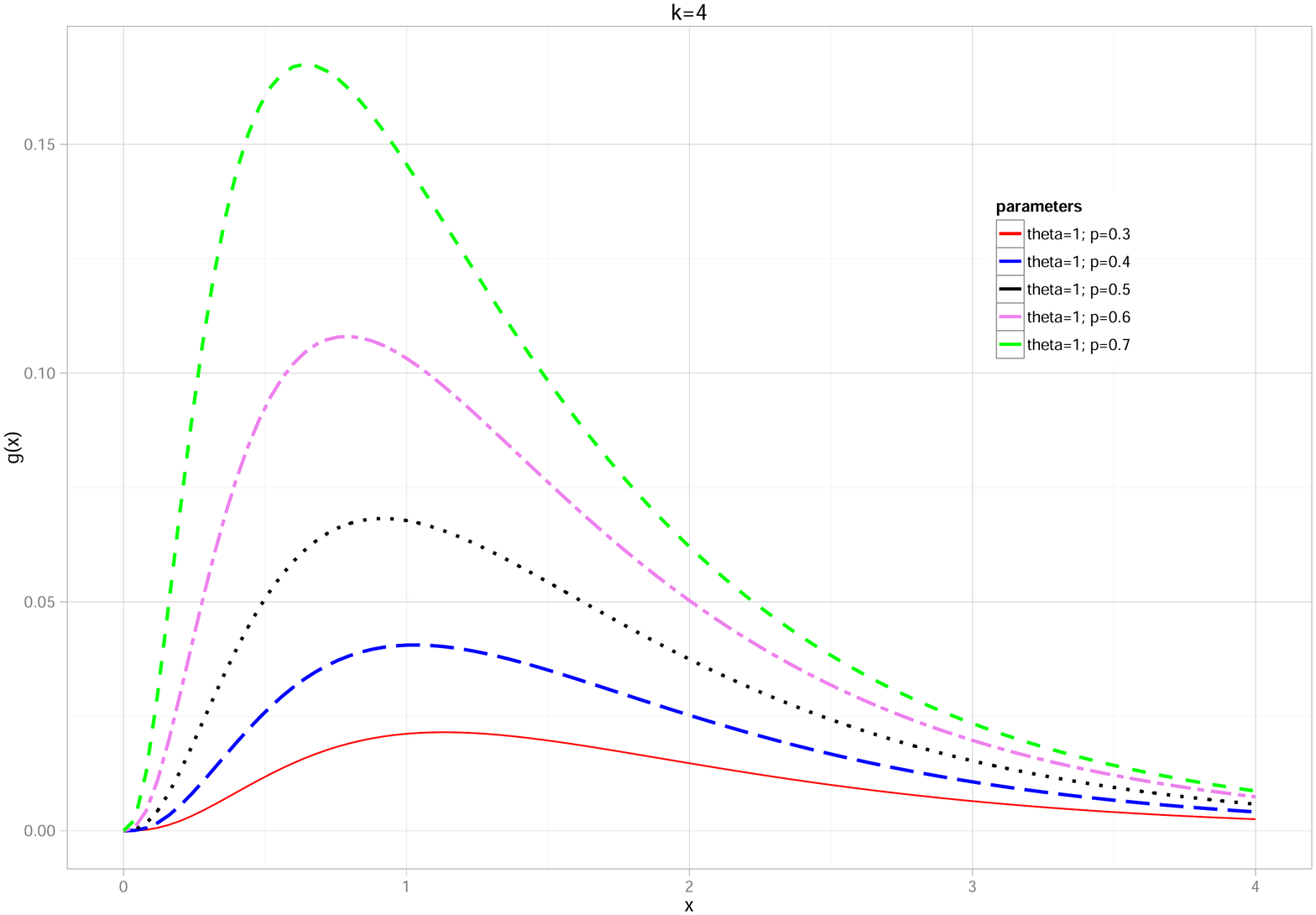}} \\
\end{tabular}
\caption{Probability functions of the EGTL distribution for $k=1,2,3,4$
\label{graph:density}}
\end{center}
\end{figure}

For $k=1$, this pdf decreases strictly in $x$ and tends to zero as $x\rightarrow \infty$. The modal value of the density of the EL distribution, at $x=0$, is given by $\frac{\theta p}{-(1-p) \ln (1-p)}$ and hence, its median is $x_\text{median}=-\frac{1}{\theta}\ln\big(\frac{1-\sqrt{1-p}}{p}\big)$. The EL distribution tends to an exponential distribution with rate parameter $\theta$, as $p\rightarrow 1$. The function is concave upward on $[0,\infty)$. The graphs of the density resemble those of the exponential and Pareto II distributions (see, Figure \ref{graph:density}).

Also, the cumulative distribution function (cdf) of $x$ corresponding to the pdf in equation (\ref{eq:marginal-x-dist}) is given by:

\begin{equation}\label{eq:cdf-x}
    G_{k}(x/p,\theta, k)=\frac{1}{A(p,k)}\sum_{r=0}^{k-1} {k-1\choose r} (-1)^{r} (1-p)^{r}
    I_{(x,r)}
\end{equation}

where,
\begin{equation*}
    I_{(x,r)}=\int_{1-p}^{1-pe^{-\theta x}}t^{-(r+1)}dt
\end{equation*}

Then, the final cdf can be reduced to:

\begin{equation} \label{eq:cdf-reduced}
\begin{split}
 G_{k}(x/p,\theta, k) & = \frac{1}{A(p,k)}\sum_{j=k}^{\infty}\frac{(py)^{j}}{j}\\
                      & = \frac{1}{A(p,k)}\Big[-\ln(1-py)-\psi(k)\sum_{j=1}^{k-1}\frac{(py)^{j}}{j}\Big]\\
               &= \frac{\ln(1-py)+\psi(k)\sum_{j=1}^{k-1}\frac{(py)^{j}}{j}}{\ln(1-p)+\psi(k)\sum_{j=1}^{k-1}\frac{p^{j}}{j}}
\end{split}
\end{equation}
where,
\[
y=\frac{1-e^{-\theta x}}{1-pe^{-\theta x}}
\]

\subsection{Moment generating function and $r^{th}$ Moment}
Suppose $x$ has the pdf in equation (\ref{eq:marginal-x-dist}), then the moment generating function (mgf) is given by:

\begin{equation}\label{eq:generating-f}
    E(e^{t x})= \frac{p^{k}}{A(p,k)}\sum_{i=0}^{\infty} {k-1+i\choose i} p^{i} \beta(i-\frac{t}{\theta}+1,k)
\end{equation}

where,
\begin{equation*}
    \beta(a,b)=\int_{0}^{1} t^{a-1} (1-t)^{b-1} d t
\end{equation*}

and hence, we can write the mgf as:

\begin{equation}\label{eq:generating-f}
    E(e^{t x})= \frac{p^{k}}{A(p,k)}\sum_{i=0}^{\infty} \sum_{j=0}^{k-1}{k-1+i\choose i} {k-1\choose j} p^{i} (-1)^{j} \frac{1}{i+j-\frac{t}{\theta}+1}
\end{equation}

\begin{eqnarray*}
  E(x) &=& \frac{\partial}{\partial t} E(e^{t x})/_{t=0}\\
    &=& \frac{p^{k}}{\theta A(p,k)}\sum_{i=0}^{\infty} \sum_{j=0}^{k-1}{k-1+i\choose i} {k-1 \choose j}p^{i} (-1)^{j} \frac{1}{(i+j-\frac{t}{\theta}+1)^2}/_{t=0} \\
    &=& \frac{p^{k}}{\theta A(p,k)}\sum_{i=0}^{\infty} \sum_{j=0}^{k-1}{k-1+i \choose i} {k-1 \choose j}p^{i} (-1)^{j} \frac{1}{(i+j+1)^{2}}
\end{eqnarray*}

The $r^{th}$ moment is given by:

\begin{equation}\label{eq:rth-moment}
    E(x^{r})= \frac{\Gamma(r+1)}{\theta^{r}}\frac{p^{k}}{A(p,k)}\sum_{i=0}^{\infty} \sum_{j=0}^{k-1}{k-1+i\choose i} {k-1\choose j} p^{i} (-1)^{j} \frac{1}{(i+j+1)^{r+1}}
\end{equation}

\subsection{Reliability and failure rate functions}\label{sec:reliability}

It is well known that the reliability (survival) function is the probability of being alive just before duration $x$, given by $S(x) = Pr(X > x) = 1 - G(x) = \int_{x}^{\infty}f(t)d t$   which is the probability that the event of interest has not occurred by duration $x$. So, the reliability $S(x)$ is the probability that a system will be successful in the interval from time $0$ to time $x$, where $X$ is a random variable denoting the time-to-failure or failure time. One may refer to the literature on reliability theory \citep{Barlow1975,Barlow1981,Basu1988}.
The survival function, corresponding to the pdf in equation (\ref{eq:marginal-x-dist}), is given by equation (\ref{eq:survivior}). Table (\ref{tab:survivior}) presents the reliability function for some special cases.

\begin{equation} \label{eq:survivior}
 S_{k}(x/p,\theta, k)  =
               1 - \frac{\ln\big(1-p\frac{1-e^{-\theta x}}{1-pe^{-\theta x}}\big)+\psi(k)\sum_{j=1}^{k-1}\frac{1}{j}
               \big(p\frac{1-e^{-\theta x}}{1-pe^{-\theta x}}\big)^{j}}{\ln(1-p)+\psi(k)\sum_{j=1}^{k-1}\frac{p^{j}}{j}}
\end{equation}

\begin{table}[htp]
\begin{center}
\small{
\caption{Reliability function for some special cases \label{tab:survivior}}
\begin{tabular}{lcc}
  \hline
  \textbf{order statistic} & $\textbf{k}$ & $\textbf{S(x)}$\\
  \hline
  &&\\
  first   & $k=1$ & $\frac{\ln (1-pe^{-\theta x})}{\ln (1-p)}$ \\
  &&\\
  second  & $k=2$ & $\frac{\ln (1-pe^{-\theta x})-p\big[\frac{1-e^{-\theta x}}{1-pe^{-\theta x}}-1\big]}{\ln (1-p)+p}$ \\
  &&\\
  third   & $k=3$ & $\frac{\ln (1-pe^{-\theta x})+p+\frac{p^{2}}{2}-p\frac{1-e^{-\theta x}}{1-pe^{-\theta x}}-
  \frac{p^{2}}{2}\big(\frac{1-e^{-\theta x}}{1-pe^{-\theta x}}\big)^{2}}
  {\ln (1-p)+p+\frac{p^{2}}{2}}$ \\
  &&\\
  \hline
\end{tabular}}
\end{center}
\end{table}

The failure rate, known as hazard rate function $h(x)$, is the instantaneous rate of occurrence of the event of interest at duration $x$ (i.e. the rate of event occurrence per unit of time). Mathematically, it is equal to the pdf of events at $x$, divided by the probability of surviving to that duration without experiencing the event. Thus, we define the failure rate function as in \cite{Barlow1963} by $h(x)=g(x)/S(x)$.
The hazard function for some special cases is given in table (\ref{tab:hazard}).

\begin{table}[htp]
\begin{center}
\small{
\caption{Failure rate function for some special cases \label{tab:hazard}}
\begin{tabular}{lcc}
  \hline
  \textbf{order statistic} & $\textbf{k}$ & $\textbf{h(x)}$\\
  \hline
    &&\\
  first   & $k=1$ & $\frac{-p \theta e^{-\theta x}}{(1-pe^{-\theta x}) \ln(1-pe^{-\theta x})}$ \\
  &&\\
  second  & $k=2$ & $\frac{-p^{2}\theta e^{-\theta x}(1-e^{-\theta x})}{(1-pe^{-\theta x})\bigg[\ln (1-pe^{-\theta x})-p\big[\frac{1-e^{-\theta x}}{1-pe^{-\theta x}}-1\big]\bigg]}$ \\
  &&\\
  third   &  $k=3$ & $\frac{-p^{3}\theta e^{-\theta x}(1-e^{-\theta x})^{2}}{(1-pe^{-\theta x})^{3}\bigg[\ln (1-pe^{-\theta x})+p+\frac{p^{2}}{2}-p\frac{1-e^{-\theta x}}{1-pe^{-\theta x}}-
  \frac{p^{2}}{2}\big(\frac{1-e^{-\theta x}}{1-pe^{-\theta x}}\big)^{2}\bigg]}$ \\
  &&\\
  \hline
\end{tabular}}
\end{center}
\end{table}

\begin{figure}[htp]
\begin{center}
\begin{tabular}[c]{cc}
\fbox{\includegraphics[scale=0.23]{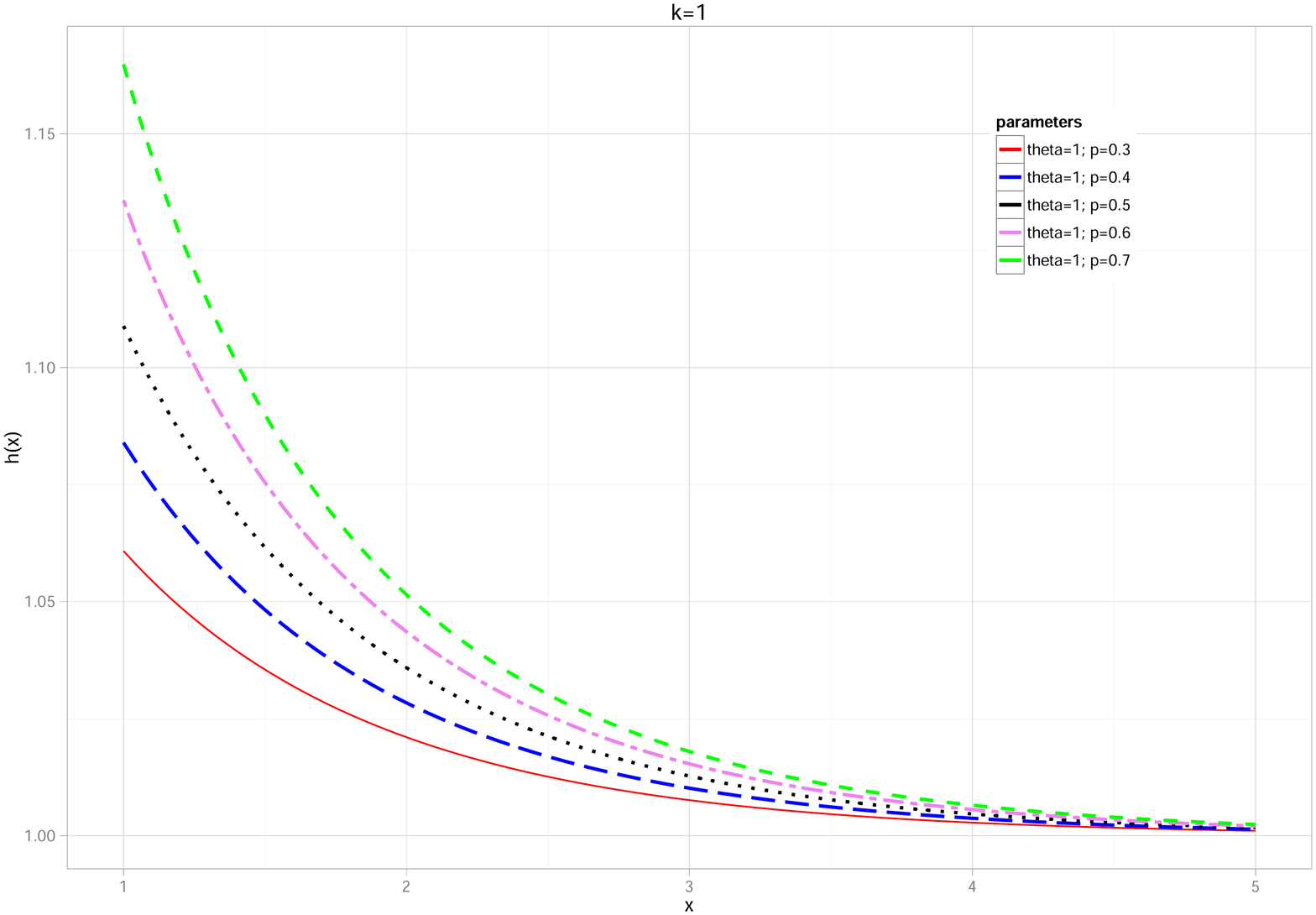}} &
\fbox{\includegraphics[scale=0.23]{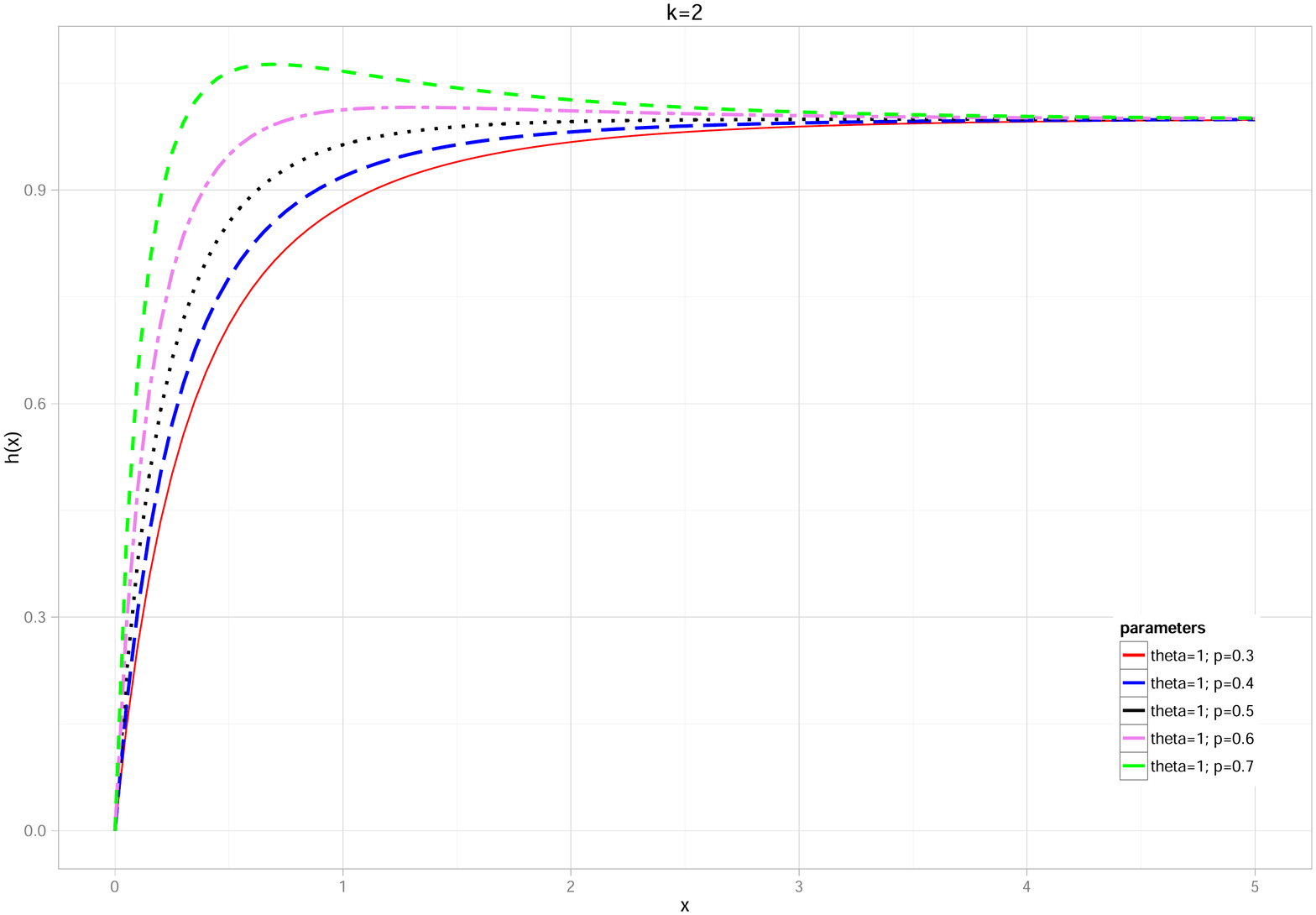}} \\
\fbox{\includegraphics[scale=0.23]{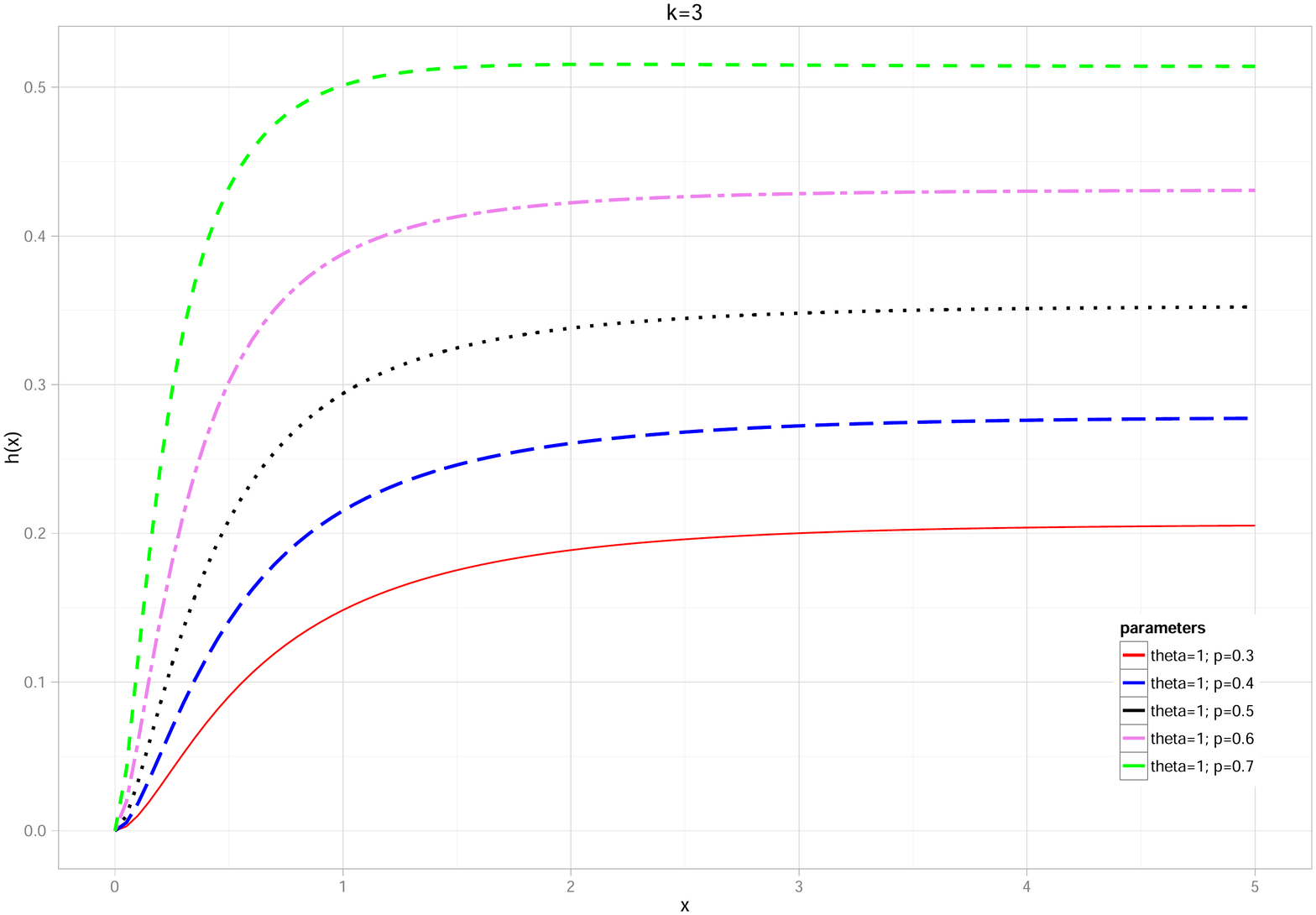}} &
\fbox{\includegraphics[scale=0.23]{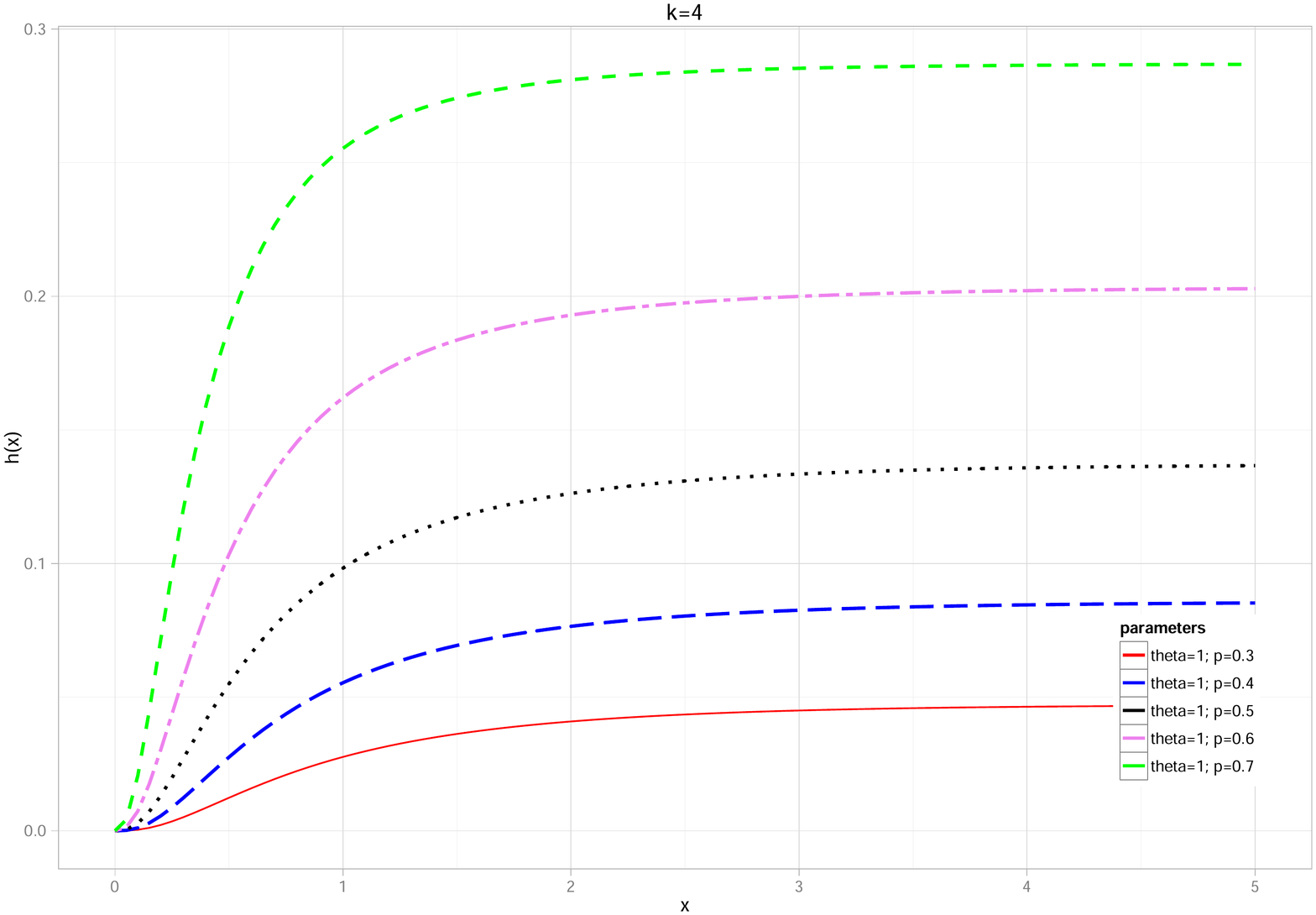}} \\
\end{tabular}
\caption{Hazard functions of the EGTL distribution for $k=1,2,3,4$
\label{graph:hazard-rate}}
\end{center}
\end{figure}

The failure rate function is analytically related to the  failure's probability distribution. It leads to the examination of the increasing (IFR) or decreasing failure rate (DFR) properties of life-length distributions. $G$ is an IFR distribution, if $h(x)$ increases for all $X$ such that $G(X)< 1$.  The motivation of the EGTL lifetime distribution is the realistic features of the hazard rate in many real-life physical and non-physical systems, which is not monotonically increasing, decreasing or constant hazard rate.
If $k=1$, the hazard rate function is decreasing following \cite{Tahmasbi2008}. In fact, if $x \rightarrow 0$ then $h(x/ p,\theta,k)=\frac{-p \theta}{(1-p) \ln (1-p)}$ and if $x \rightarrow \infty$ then  $h(x/p,\theta ,k) \rightarrow  \theta$.
For $k>1$, there is an increasing failure rate. Indeed, if $x \rightarrow 0$ then $h(x/p,\theta ,k) \rightarrow  0$. If $x \rightarrow \infty$ then  $h(x/p,\theta ,k) >  0$ (see Figure \ref{graph:hazard-rate}).

\subsection{Random number generation}\label{sec:random-generation}

We can generate a random variable from the cdf of $x$ in equation (\ref{eq:cdf-reduced}) using the following steps:

\begin{itemize}
  \item Generate a random variable $U$ from the standard uniform distribution.
  \item Solve the non linear equation in $y$:

      \begin{equation} \label{eq:generate-y}
        U = \frac{\ln(1-py)+\psi(k)\sum_{j=1}^{k-1}\frac{(py)^{j}}{j}}{\ln(1-p)+\psi(k)\sum_{j=1}^{k-1}\frac{p^{j}}{j}}
      \end{equation}

     where,

      \begin{equation*}\label{eq:psi-k}
        \psi(k) =
          \begin{cases}
           0 & \mbox{if } k=1 \\
           1 & \mbox{if } k=2, 3, ...,z
          \end{cases}
      \end{equation*}

  \item Calculate the values of $X$ such as:
       \begin{equation} \label{eq:generate-x}
       X=-\frac{1}{\theta}\ln\bigg(\frac{1-y}{1-py}\bigg)
       \end{equation}
\end{itemize}

where $X$ is EGTL random variable with parameters $\theta$ and $p$. Note that for the special case $k=1$, we generate $X$ directly from the following equation:

\begin{equation}\label{generate-x2}
    X=-\frac{1}{\theta}\ln\bigg(\frac{1-(1-p)^{1-u}}{p}\bigg)
\end{equation}

\section{Estimation of the parameters}

In this section, we will determine the estimates of the parameters $p$ and $\theta$ for the EGTL new family of distributions. There are many methods available for estimating the parameters of interest. We present here the three most popular methods: Maximum likelihood, method of moments and Bayesian estimations.

\subsection{Maximum Likelihood estimation}

Let ($X_{1},X_{2}, \dots,X_{n}$) be a random sample from the EGTL distribution. The log-likelihood function given the observed values, $x_{obs}=(x_{1},x_{2}, \dots,x_{n})$, is:

\begin{multline}
    \ln L= n\ln(\theta)+nk\ln(p)-\theta \sum_{i=1}^{n}x_{i}+(k-1)\sum_{i=1}^{n}\ln(1-e^{-\theta x_{i}})
    \\ -n\ln A(p,k)-k\sum_{i=1}^{n}\ln(1-pe^{-\theta x_{i}})
\end{multline}

We subsequently derive the associated gradients:

\[
\frac{\partial \ln L}{\partial p}
  = \frac{np^{k-1}}{(1-p)\big[\ln(1-p)+\psi(k)\sum_{j=1}^{k-1}\frac{p^{j}}{j}\big]}+\frac{nk}{p}+\sum_{i=1}^{n}\frac{1}{e^{\theta x_{i}}-p}
\]

\[
\frac{\partial \ln L}{\partial \theta}
=
\frac{n}{\theta}-\sum_{i=1}^{n}x_{i}+(k-1)\sum_{i=1}^{n}\frac{x_{i}}{e^{\theta x_{i}}-1}-kp\sum_{i=1}^{n}\frac{x_{i}}{e^{\theta x_{i}}-p}
\]

We need the Fisher information matrix for interval estimation and tests of hypotheses on the parameters.
It can be expressed in terms of the second derivatives of the log-likelihood function:

\[
\mathcal{I}(\widehat{p},\widehat{\theta})
  = -
\left.
\left( \begin{array}{cc}
  E\Big(\tfrac{\partial^2 \ln L}{\partial p^2}\Big)
  &  E\Big(\tfrac{\partial^2 \ln L}{\partial p \partial \theta}\Big) \\
  E\Big(\tfrac{\partial^2 \ln L}{\partial \theta \partial p}\Big)
  &  E\Big(\tfrac{\partial^2 \ln L}{\partial \theta^2}\Big) \\
\end{array} \right)
\right|_{\theta = (\widehat{p},\widehat{\theta})}
\]

The maximum likelihood estimates (MLEs) $\widehat{p}$ and $\widehat{\theta}$  of  the EGTL parameters $p$ and $\theta$, respectively, can be found analytically using the iterative EM algorithm to handle the incomplete data problems \citep{Dempster1977,Krishnan1997}. The iterative method consists in repeatedly updating the parameter estimates by replacing the "missing data" with the new estimated values. The standard method used to determine the MLEs is the Newton-Raphson algorithm that requires second derivatives of the log-likelihood function for all iterations. The main drawback of the EM algorithm is its rather slow convergence, compared to the Newton-Raphson method, when the "missing data" contain a relatively large amount of information \citep{Little1983}. Recently, several researchers have used the EM method such as \cite{Adamidis1998}, \cite{Adamidis2005}, \cite{Karlis2003}, \cite{Ng2002} and others. Newton-Raphson is required for the M-step of the EM algorithm.
To start the algorithm, a hypothetical distribution of complete-data is defined with the pdf in equation (\ref{eq:joint-dist}) and then, we drive the conditional mass function as:

\begin{equation}\label{eq:prob-z}
   p(z/x,p,\theta)=\frac{\Gamma(z)}{\Gamma(k)\Gamma(z-k+1)}p^{z-k} e^{-\theta (z-k)x}(1-pe^{-\theta x})^{k}
\end{equation}

\textbf{E-step:}

\begin{equation}\label{eq:E-step}
   E(z/x,p,\theta)=\frac{k}{1-pe^{-\theta x}}
\end{equation}

\textbf{M-step:}

\begin{equation}\label{eq:M-step1}
\widehat{p}^{(r+1)}
  =
  \frac{-\big(1-p^{(r+1)}\big)\Big[\ln\big(1-p^{(r+1)}\big)+\psi(k)\sum_{j=1}^{k-1}\frac{\big(p^{(r+1)}\big)^{j}}{j}\Big]}{n\big(p^{(r+1)}\big)^{k-1}}\sum_{i=1}^{n}\frac{k}{1-p^{(r)}e^{-\theta^{(r)} x_{i}}}
\end{equation}

\begin{equation}\label{eq:M-step2}
\widehat{\theta}^{(r+1)}  =
n\bigg[\sum_{i=1}^{n}\frac{kx_{i}}{1-p^{(r)}e^{-\theta^{(r)} x_{i}}}-(k-1)\sum_{i=1}^{n}\frac{x_{i}}{1-e^{-\theta^{(r+1)} x_{i}}}\bigg]^{-1}
\end{equation}

\subsection{Method of moments estimation}

The method of moments involves equating theoretical with sample moments. The estimate of $r^{th}$ moment is
$\widehat{\mu}_{r}=\frac{1}{n}\sum_{i=1}^{n}x_{i}^{r}$. For the EGTL distribution, the $r^{th}$ moment is given by equation (\ref{eq:rth-moment}). The corresponding first and second moments are given by:

\begin{equation}\label{eq:first-moment}
    \frac{1}{n}\sum_{i=1}^{n}x_{i}=E(x)= \frac{1}{\theta}\frac{p^{k}}{A(p,k)}\sum_{i=0}^{\infty} \sum_{j=0}^{k-1}{k-1+i\choose i} {k-1\choose j} p^{i} (-1)^{j} \frac{1}{(i+j+1)^{2}}
\end{equation}

\begin{equation}\label{eq:second-moment}
    \frac{1}{n}\sum_{i=1}^{n}x_{i}^{2}= E(x^{2})= \frac{2}{\theta^{2}}\frac{p^{k}}{A(p,k)}\sum_{i=0}^{\infty} \sum_{j=0}^{k-1}{k-1+i\choose i} {k-1 \choose j}p^{i} (-1)^{j} \frac{1}{(i+j+1)^{3}}
\end{equation}

From equation (\ref{eq:first-moment}) we obtain:

\begin{equation}\label{eq:p-moment}
    \theta = \frac{p^{k}}{\overline{x}A(p,k)}\sum_{i=0}^{\infty} \sum_{j=0}^{k-1}{k-1+i\choose i} {k-1\choose j}p^{i} (-1)^{j} \frac{1}{(i+j+1)^{2}}
\end{equation}

and then, we should solve the following equation in $p$:

\begin{equation}\label{eq:theta-moment}
   \frac{1}{n}\sum_{i=1}^{n}x_{i}^{2}-\frac{2\overline{x}^{2}A(p,k)\sum_{i=0}^{\infty} \sum_{j=0}^{k-1}{k-1+i\choose i} {k-1 \choose j}p^{i} (-1)^{j} \frac{1}{(i+j+1)^{3}}} {p^{k}\big[\sum_{i=0}^{\infty} \sum_{j=0}^{k-1}{k-1+i\choose i} {k-1 \choose j}p^{i} (-1)^{j} \frac{1}{(i+j+1)^{2}}\big]^2}=0
\end{equation}

Thereafter, we determine $\widehat{\theta}$ by replacing $p$ with its estimated value, $\widehat{p}$, in the equation (\ref{eq:p-moment}).

\subsection{Bayesian estimation}

In the Bayesian approach inferences are expressed in a posterior distribution for the
parameters which is, according to Bayes' theorem, given in terms of the likelihood and a prior density function by:

\begin{equation}\label{eq:bayes-posterior}
    P_{k}(p,\theta/x_{1}, x_{2}, ..., x_{n})=\frac{g_{k}(x/p,\theta).\pi_{k}(p,\theta)}{g_{k}(x)}
\end{equation}

where, $\pi_{k}(p,\theta)$ is a prior probability distribution function and $g_{k}(x/p,\theta)$
is the likelihood of observations $(x_{1}, x_{2}, ..., x_{n})$. Note that $g_{k}(x)$ is the normalizing constant for the function $\pi_{k}(p,\theta)g_{k}(x/p,\theta)$ given by:

\begin{equation}\label{eq:bayes-normalizing}
    \int\int \pi_{k}(p,\theta)g_{k}(x/p,\theta) dp d\theta
\end{equation}

We should first specify our initial beliefs or other sorts of knowledge on the prior distribution $\pi_{k}(p,\theta)$. Here, we suppose that the standard uniform distribution on the interval $[0, 1]$ is a prior distribution for the parameter $p$ and gamma, $G(a,b)$, is a prior distribution for the parameter $\theta$, where $a$ is a shape parameter and $b$ is a scale parameter. The prior probability function is then equal to:

\begin{equation}\label{eq:bayes-prior}
    \pi_{k}(p,\theta)=\pi_{1,k}(p)\pi_{2,k}(\theta)
\end{equation}

where, $\pi_{1,k}(p)=1$ and $\pi_{2,k}(\theta)=\frac{b^a}{\Gamma(a)}\theta^{a-1}e^{-b\theta}$.

Using the mean square error as a risk function, we obtain the Bayes estimates as the means of the posterior distribution:

\begin{equation}\label{eq:Bayes-p-estimate}
    \widehat{\theta} =E(\theta)= \int_{0}^{\infty}\int_{0}^{1}\theta P_{k}(p,\theta/x)d\theta dp
\end{equation}

\begin{equation}\label{eq:Bayes-theta-estimate}
    \widehat{p} =E(p)= \int_{0}^{1}\int_{0}^{\infty}p P_{k}(p,\theta/x)dp d\theta
\end{equation}

\section{Simulation study}

As an illustration of the three last methods of estimation, numerical computations have been performed using the steps presented in section \ref{sec:random-generation} for the random number generation. The numerical study was
based on $1000$ random samples of the sizes $20$, $50$ and $100$  from the EGTL distribution for each of the $3$ values of $\lambda=(p,\theta)$ and the three cases $k=\{1, 2, 3\}$.  We have considered the initial values $(0.5 ; 0.5)$, $(0.7 ; 1.5)$ and $(0.3 ; 2)$. For this purpose, we have used the program
Mathcad 14.0.
After determining the parameter estimates $\widehat{\lambda}=(\widehat{p},\widehat{\theta})$ we compute the biases, the variances and the mean square errors (MSEs),  where
$MSE(\widehat{\lambda})= E(\widehat{\lambda}-\lambda)^{2}= Bias^{2}(\widehat{\lambda}) + var(\widehat{\lambda})$
and $Bias(\widehat{\lambda}) = E(\widehat{\lambda})-\lambda$.
An estimator $\widehat{\lambda}$ is said to be efficient if its mean square error (MSE) is minimum
among all competitors. In fact, $\widehat{\lambda}_{1}$ is more efficient than $\widehat{\lambda}_{2}$ if $MSE(\widehat{\lambda}_{1}) < MSE(\widehat{\lambda}_{2})$.

Table \ref{tab:simulation} reports the results from the simulated data where the variances and the MSEs of the parameters are given. The results show that, for each case $k=\{1, 2, 3\}$, the variances and the MSEs decrease when the sample size increases.
We see that the values from the Bayesian method are generally lower than those obtained using the ML approach.

\begin{sidewaystable}[htp]
  \centering
  \caption{Results from the numerical computation}
  \scriptsize{
    \begin{tabular}{lllcccccccccccc}
          &       &       & \multicolumn{4}{c}{\textbf{Maximum likelihood }} & \multicolumn{4}{c}{\textbf{Method of Moments}} & \multicolumn{4}{c}{\textbf{Bayesian methods}}\\
\cline{4-15}    \textbf{n} & \textbf{k} & \textbf{$(p,\theta)$} & \textbf{$var(\widehat{p})$} & \textbf{$var(\widehat{\theta})$} & \textbf{$MSE(\widehat{p})$} & \textbf{$MSE(\widehat{\theta})$} & \textbf{$var(\widehat{p})$} & \textbf{$var(\widehat{\theta})$} & \textbf{$MSE(\widehat{p})$} & \textbf{$MSE(\widehat{\theta})$} & \textbf{$var(\widehat{p})$} & \textbf{$var(\widehat{\theta})$} & \textbf{$MSE(\widehat{p})$} & \textbf{$MSE(\widehat{\theta})$} \\
    \hline
    \textbf{n=20} &       &       &       &       &       &       &       &       &       &       &       &       &       &  \\
          & \textbf{1} & (0.5 ; 0.5) & 0.2091 & 0.2008 & 0.2136 & 0.2043 & 0.2775 & 0.2241 & 0.2849 & 0.2351 & 0.1737 & 0.1414 & 0.1835 & 0.1586 \\
          &       & (0.7 ; 1.5) & 0.2779 & 0.6980 & 0.2861 & 0.7145 & 0.3066 & 1.0227 & 0.3078 & 1.0361 & 0.2208 & 0.5794 & 0.2287 & 0.6104 \\
          &       & (0.3 ; 2) & 0.4250 & 1.0930 & 0.4680 & 1.1160 & 0.6890 & 1.7260 & 0.7270 & 1.8930 & 0.2110 & 0.5010 & 0.2420 & 0.5870 \\
          &       &       &       &       &       &       &       &       &       &       &       &       &       &  \\
          & \textbf{2} & (0.5 ; 0.5) & 0.2725 & 0.2013 & 0.2773 & 0.2027 & 0.3087 & 0.2327 & 0.3347 & 0.2398 & 0.1934 & 0.1227 & 0.1957 & 0.1488 \\
          &       & (0.7 ; 1.5) & 0.1105 & 1.1326 & 0.1114 & 1.1380 & 0.1318 & 1.2050 & 0.1409 & 1.2078 & 0.0742 & 0.9730 & 0.0796 & 1.0117 \\
          &       & (0.3 ; 2) & 0.2700 & 0.6230 & 0.3660 & 0.6520 & 0.5750 & 1.0000 & 0.6810 & 1.1900 & 0.1000 & 0.3690 & 0.1230 & 0.4320 \\
          &       &       &       &       &       &       &       &       &       &       &       &       &       &  \\
          & \textbf{3} & (0.5 ; 0.5) & 0.1407 & 0.1164 & 0.1193 & 0.1536 & 0.2186 & 0.1498 & 0.2245 & 0.1515 & 0.1054 & 0.0490 & 0.1065 & 0.0533 \\
          &       & (0.7 ; 1.5) & 0.1782 & 0.6809 & 0.1938 & 0.7694 & 0.2247 & 0.9003 & 0.2269 & 0.9290 & 0.1653 & 0.5840 & 0.1682 & 0.5880 \\
          &       & (0.3 ; 2) & 0.4090 & 0.8630 & 0.4410 & 0.9040 & 1.2390 & 0.6680 & 1.2930 & 0.7200 & 0.1590 & 0.3380 & 0.2090 & 0.3810 \\
    \textbf{n=50} &       &       &       &       &       &       &       &       &       &       &       &       &       &  \\
          & \textbf{1} & (0.5 ; 0.5) & 0.1099 & 0.1032 & 0.1119 & 0.1049 & 0.1534 & 0.1280 & 0.1848 & 0.1308 & 0.0681 & 0.0493 & 0.0723 & 0.0506 \\
          &       & (0.7 ; 1.5) & 0.1221 & 0.5182 & 0.1257 & 0.5431 & 0.1426 & 0.6526 & 0.1532 & 0.6929 & 0.0567 & 0.3863 & 0.0625 & 0.3878 \\
          &       & (0.3 ; 2) & 0.3630 & 0.8440 & 0.3820 & 0.9340 & 0.6330 & 1.6030 & 0.6580 & 1.6480 & 0.1380 & 0.3710 & 0.1730 & 0.4240 \\
          &       &       &       &       &       &       &       &       &       &       &       &       &       &  \\
          & \textbf{2} & (0.5 ; 0.5) & 0.1687 & 0.1097 & 0.1915 & 0.1239 & 0.2252 & 0.1422 & 0.2257 & 0.1544 & 0.1014 & 0.0443 & 0.1119 & 0.0456 \\
          &       & (0.7 ; 1.5) & 0.1051 & 1.0626 & 0.1072 & 1.0985 & 0.1287 & 1.1868 & 0.1381 & 1.1923 & 0.0691 & 0.9795 & 0.0735 & 0.9807 \\
          &       & (0.3 ; 2 & 0.2350 & 0.5620 & 0.2510 & 0.5770 & 0.5330 & 0.8980 & 0.5650 & 0.9390 & 0.0340 & 0.2670 & 0.0740 & 0.3050 \\
          &       &       &       &       &       &       &       &       &       &       &       &       &       &  \\
          & \textbf{3} & (0.5 ; 0.5) & 0.0941 & 0.0782 & 0.0947 & 0.0896 & 0.1782 & 0.1134 & 0.1810 & 0.1167 & 0.0624 & 0.0185 & 0.0686 & 0.0224 \\
          &       & (0.7 ; 1.5) & 0.0606 & 0.5371 & 0.0608 & 0.5459 & 0.1005 & 0.6989 & 0.1067 & 0.7072 & 0.0433 & 0.4644 & 0.0438 & 0.4685 \\
          &       & (0.3 ; 2 & 0.3360 & 0.7130 & 0.3790 & 0.8170 & 0.5770 & 1.1460 & 0.6340 & 1.1800 & 0.1310 & 0.2140 & 0.1430 & 0.2760 \\
    \textbf{n=100} &       &       &       &       &       &       &       &       &       &       &       &       &       &  \\
          & \textbf{1} & (0.5 ; 0.5) & 0.0652 & 0.0398 & 0.0663 & 0.0592 & 0.1049 & 0.0827 & 0.1098 & 0.0879 & 0.0185 & 0.0022 & 0.0327 & 0.0082 \\
          &       & (0.7 ; 1.5) & 0.0531 & 0.4206 & 0.0605 & 0.4663 & 0.0870 & 0.5859 & 0.0901 & 0.5876 & 0.0013 & 0.1941 & 0.0052 & 0.2201 \\
          &       & (0.3 ; 2) & 0.2588 & 0.6168 & 0.2828 & 0.7518 & 0.4918 & 1.2348 & 0.5318 & 1.3888 & 0.0118 & 0.0058 & 0.0618 & 0.1188 \\
          &       &       &       &       &       &       &       &       &       &       &       &       &       &  \\
          & \textbf{2} & (0.5 ; 0.5) & 0.1059 & 0.0733 & 0.1140 & 0.0804 & 0.1897 & 0.1097 & 0.1971 & 0.1162 & 0.0479 & 0.0051 & 0.0574 & 0.0185 \\
          &       & (0.7 ; 1.5) & 0.0522 & 0.9943 & 0.0614 & 1.0109 & 0.0849 & 1.1288 & 0.0889 & 1.1395 & 0.0002 & 0.8165 & 0.0039 & 0.8951 \\
          &       & (0.3 ; 2) & 0.1739 & 0.4719 & 0.2009 & 0.5239 & 0.4229 & 0.7759 & 0.4979 & 0.8269 & 0.0049 & 0.0149 & 0.0219 & 0.0859 \\
          &       &       &       &       &       &       &       &       &       &       &       &       &       &  \\
          & \textbf{3} & (0.5 ; 0.5) & 0.0749 & 0.0560 & 0.0814 & 0.0613 & 0.1255 & 0.0888 & 0.1471 & 0.0901 & 0.0247 & 0.0005 & 0.0505 & 0.0059 \\
          &       & (0.7 ; 1.5) & 0.0174 & 0.4517 & 0.0199 & 0.4680 & 0.0543 & 0.6477 & 0.0586 & 0.6599 & 0.0006 & 0.3827 & 0.0011 & 0.4055 \\
          &       & (0.3 ; 2) & 0.2183 & 0.4733 & 0.2873 & 0.5973 & 0.4423 & 1.0013 & 0.4923 & 1.0643 & 0.0043 & 0.0073 & 0.0663 & 0.0323 \\
    \hline
    \end{tabular}}
  \label{tab:simulation}
\end{sidewaystable}

\section{Application examples}

In this section, we fit the EGTL distribution to two real data sets using the MLEs. The first set (table \ref{tab:data-Barlow-Campo}) consists of "$107$ failure times for right rear brakes on D9G-66A caterpillar tractors", reproduced from \cite{Barlow1975a} and used also by \cite{Chang1993}. These data are used in many applications of reliability \citep{Adamidis2005,Tsokos2012,Shahsanaei2012}.
The second set of data involves $100$ observations (table \ref{tab:data-Quesenberry}) of the results from an experiment concerning "the tensile fatigue characteristics of a polyester/viscose yarn". These data were presented by \cite{Picciotto1970} to study the problem of warp breakage during weaving.  The observations were obtained on the cycles to failure of a $100$ cm yarn sample put to test under $2.3\%$ strain level. The sample is used in \cite{Quesenberry1982} as an example to illustrate selection procedure among probability distributions used in reliability. The reliability function of these two data sets belongs to the increasing failure rate class \citep{Doksum1984,Adamidis2005}. In addition to our class of distributions, the gamma and Weibull distributions were fitted these data sets. The respective densities of gamma and Weibull distributions are $f_{1}(x,\lambda_{1},\beta_{1})=\lambda_{1}^{\beta_{1}}x^{\beta_{1}-1}exp(-\lambda_{1}x)\Gamma(\beta_{1})^{-1}$ and $f_{2}(x,\lambda_{2},\beta_{2})=\beta_{2}\lambda_{2}^{\beta_{2}}x^{\beta_{2}-1}exp(-\lambda_{2}x)^{\beta_{2}}$.

\begin{table}[htp]
\begin{center}
\scriptsize{
\caption{"Ordered Failure Times (in hours) of 107 Right Rear Brakes on D9G-66A Caterpillar Tractors" \cite{Barlow1975a,Chang1993} \label{tab:data-Barlow-Campo}}
\begin{tabular}{lllllllllllllll}
  \hline
56&	753&	1153&	1586&	2150&	2624&	3826&	83&	763&	1154&	1599&	2156&	2675&	3995&	104\\
806&	1193&	1608&	2160&	2701&	4007&	116&	834&	1201&	1723&	2190&	2755&	4159&	244&	838\\
1253&	1769&	2210&	2877&	4300&	305&	862&	1313&	1795&	2220&	2879&	4487&	429&	897&	1329\\
1927&	2248&	2922&	5074&	452	&904&	1347&	1957&	2285&	2986&	5579&	453&	981&	1454&	2005\\
2325&	3092&	5623&	503&	1007&	1464&	2010&	2337&	3160&	6869&	552&	1008&	1490&	2016&	2351\\
3185&	7739&	614&	1049&	1491&	2022&	2437&	3191&	661&	1069&	1532&	2037&	2454&	3439&	673\\
1107&	1549&	2065&	2546&	3617&	683&	1125&	1568&	2096&	2565&	3685&	685&	1141&	1574&	2139\\
2584&	3756&&&&&&&&&&&&&\\													
 \hline
\end{tabular}}
\end{center}
\end{table}

\begin{table}[htp]
\begin{center}
\scriptsize{
\caption{"Results of Model Selection Program on Yarn Data" \citep{Quesenberry1982} \label{tab:data-Quesenberry}}
\begin{tabular}{lllllllllllllll}
  \hline
86&	146&	251&	653	&98&	249&	400&	292&	131&	169&175&	176&	76&	264&	15\\	
364&	195&	262&	88&	264& 157&	220&	42&	321&	180&	198&	38&	20&	61&	121\\
282&	224&	149	&180&	325&	250&	196&	90&	229&	166&38&	337&	65&	151&	341\\
40&	40&	135&	597&	246& 211&	180&	93&	315&	353&	571&	124&	279&	81&	186\\
497&	182&	423&	185&	229&	400&	338&	290&	398&	71& 246&	185&	188&	568& 55\\	
55&	61&	244&	20&	284& 393&	396&	203&	829&	239&	286&	194&	277&	143&	198\\
264&	105	&203&	124&	137&	135&	350&	193	&188&	236&&&&&\\
 \hline
\end{tabular}}
\end{center}
\end{table}

Table \ref{tab:fitted} shows the fitted parameters, the calculated values of Kolmogorov-Smirnov (K-S) and their respective \emph{p-values} for the two sets of data. It should be noted that the K-S test compares
an empirical distribution with a known (not estimated) one. It is used to decide if a sample comes from a population with a specific distribution ($H_{0}$: the data follow a specified distribution). We estimate some special cases ($k=1, 2, 3, 4$) of the EGTL family of distributions at $5\%$ significant level. The \emph{p-values} are only significant for the case $k=1$ for the \cite{Barlow1975a} and \cite{Quesenberry1982} data sets. In fact, the data exhibit increasing failure rates but, the EGTL distribution is a decreasing failure rate if $k=1$ (see figure \ref{graph:hazard-rate}). The new lifetime distribution provides good fit to the data sets. The K-S test shows that the EGTL distribution is an attractive alternative to the popular gamma and Weibull distributions. It generalizes the reliability lifetime distributions to any $k^{th}$ order statistics.
Indeed, as shown in section \ref{sec:reliability}, If $k=1$, the hazard rate function is decreasing following \cite{Tahmasbi2008} and there is an increasing hazard rate for $k>1$.

\begin{table}[htp]
\begin{center}
\small{
\caption{The Goodness of Fit for some Special Cases \label{tab:fitted}}
\begin{tabular}{|l|c|c|c|c|}
  \hline
  Distributions & $\widehat{p}$ & $\widehat{\theta}$ & K-S value & p-value \\
  \hline
\textbf{\cite{Barlow1975a} data set ($n=107$):} &&&&\\
	First order (k=1)   &$0.0500$   &$5.00$  $10^{-6}$    & $0.9611$  & $0.0000$ \\
	Second order (k=2)  &$0.0232$   &$7.32$  $10^{-4}$    &$0.0639$   &$0.7746$  \\
	Third order (k=3)   &$0.8811$   &$4.38$  $10^{-4}$    &$0.1106$   &$0.1456$  \\
	Fourth order (k=4)  &$0.4209$   &$8.84$  $10^{-4}$    &$0.0723$   &$0.6305$  \\
    Gamma    & \multicolumn{2}{|c|}{$(0.943;1.908)$}& $0.0680$  & $0.7343$\\
    Weibull  & \multicolumn{2}{|c|}{$(0.447;1.486)$}& $0.0490$  & $0.9999$\\
 \hline
 \textbf{\cite{Quesenberry1982} data set ($n=100$):} &&&&\\
	First order (k=1)   &$0.1901$   &$4.22$  $10^{-3}$    & $0.1955$  & $0.0009$ \\
	Second order (k=2)  &$0.0248$   &$6.65$  $10^{-3}$    &$0.1078$   &$0.1952$  \\
	Third order (k=3)   &$0.2127$   &$7.66$  $10^{-3}$    &$0.0879$   &$0.4218$  \\
	Fourth order (k=4)  &$0.1031$   &$9.10$  $10^{-3}$    &$0.0786$   &$0.5657$  \\
    Gamma    & \multicolumn{2}{|c|}{$(1.008;2.239)$}& $0.0950$  & $0.3118$\\
    Weibull  & \multicolumn{2}{|c|}{$(0.403;1.604)$}& $0.0760$  & $0.6080$\\
 \hline
\end{tabular}}
\end{center}
\end{table}

\section{Conclusion}
We define a new two-parameter lifetime distribution so-called EGTL distribution. Our procedure generalizes the EL distribution proposed by \cite{Tahmasbi2008}. We derive some mathematical properties and we present the plots of the pdf and the failure rate functions for some special cases. The estimation of the parameters is attained by the maximum likelihood, EM algorithm, the method of moments and the Bayesian approach, with numerical computations performed as illustration of the different methods of estimation. The application study is illustrated based on two real data sets used in many applications of reliability.
We have shown that our proposed EGTL distribution is suitable for modelling the time to any failure and not only the time to the first or the last failure. It is very competitive compared with its standard counterpart's distributions.

Ordered random variables are already known for their ascending order. The paper may be extended to the concept of dual generalized ordered statistics, introduced by \cite{Burkschat2003}, that enables a common approach to the descending ordered spacings like the reverse ordered statistics and the lower record values.


\section*{Appendix}
\small

Let $T=(T_{1}, T_{2}, ..., T_{Z})$ be iid exponential r.v. with pdf given by: $f(t)=\theta e^{-\theta t}$ , for $t>0$, where $Z$  is a log-series r.v. with pmf, $P(Z=z)$, given by:

\begin{equation}\label{eq:logarithmic-appendix}
   P(Z=z) = \frac{1}{-\ln(1-p)} \frac{p^z}{z} ; z \in \{1,2,3,\dots\}
\end{equation}

From the Taylor series, for $|x|<1$ we have:
\begin{equation*}
   ln(1+x)=\sum_{j=1}^{\infty}\frac{(-1)^{j+1}}{j}x^{j}
\end{equation*}

then,
\begin{equation*}
    \ln(1-p)= - \sum_{j=1}^{\infty}\frac{p^{j}}{j}
\end{equation*}

The truncated at $k-1$ logarithmic distribution with parameter $p$ is:

\begin{equation}\label{eq:tr-logarithmic-appendix}
   P_{k-1}(Z=z) = \frac{1}{A(p,k)} \frac{p^{z}}{z} ; k = 1,2,3,\dots\ , z \mbox{ and } z=k, k+1,\dots\
\end{equation}

where,
\begin{equation}\label{eq:tr-logarithmic-A-appendix}
   A(p,k) = \sum_{j=k}^{\infty}\frac{p^{j}}{j}
          = - \ln(1-p)-\psi(k)\sum_{j=1}^{k-1}\frac{p^{j}}{j}
\end{equation}

and,
\begin{equation}\label{eq:psi-k-appendix}
   \psi(k) =
   \begin{cases}
   0 & \mbox{if } k=1 \\
   1 & \mbox{if } k=2, 3, ...,z
   \end{cases}
\end{equation}

The pdf of the $k^{th}$  order statistic is:

\begin{equation}\label{eq:order-appendix}
   f_{k}(x/z,\theta)=\frac{\theta \Gamma(z+1)}{\Gamma(k)\Gamma(z-k+1)}e^{-\theta(z-k+1)x}(1-e^{-\theta x})^{k-1} \mbox{ ; } \theta, x > 0
\end{equation}

Then, the joint distribution from eq. (\ref{eq:tr-logarithmic-appendix}) and (\ref{eq:order-appendix}) is:

\begin{equation}\label{eq:joint-dist-appendix}
\begin{split}
g_{k}(x,z/p,\theta) & =f_{k}(x/z,\theta)P_{k-1}(z/p)\\
                    & = \frac{\theta \Gamma(z+1)}{\Gamma(k)\Gamma(z-k+1)}e^{-\theta(z-k+1)x}(1-e^{-\theta x})^{k-1}\frac{1}{A(p,k)} \frac{p^{z}}{z}\\
                    & = \frac{\Gamma(z)}{\Gamma(k)\Gamma(z-k+1)}\frac{\theta p^{z} e^{-\theta (z-k+1) x} (1-e^{-\theta x})^{k-1}}{A(p,k)}
\end{split}
\end{equation}

Let, $f=\theta e^{-\theta x}$,  $F=1- e^{-\theta x}$, and  $a=pe^{-\theta x}=p(1-F)$

\begin{equation}
    g_{k}(x,z/p,\theta) = \frac{1}{A(p,k)}\frac{(z-1)!p^k}{(z-k)!(k-1)!}fF^{k-1}a^{z-k}
\end{equation}

the marginal density of $x$ is:
\begin{equation}
g_{k}(x/p,\theta) = \frac{p^{k}fF^{k-1}}{A(p,k)}\sum_{z=k}^{\infty}\frac{(z-1)!}{(z-k)!(k-1)!}a^{z-k}
\end{equation}

Let, $z-k=s$ , $z=k+s$, $k-1=z-s-1$

\begin{equation}
\begin{split}
g_{k}(x/p,\theta) &= \frac{p^{k}fF^{k-1}}{A(p,k)}\sum_{s=0}^{\infty}\frac{(s+k-1)!}{s!(k-1)!}a^{s}\\
                  &= \frac{p^{k}fF^{k-1}}{A(p,k)}\sum_{s=0}^{\infty}{s+k-1\choose s}a^{s}\\
                  &= \frac{p^{k}fF^{k-1}}{A(p,k)}\frac{1}{(1-a)^{k}}\\
                  &=   \frac{\theta p^{k} e^{-\theta x} (1-e^{-\theta x})^{k-1}}{A(p,k)(1-pe^{-\theta x})^{k}}
    \mbox{ ; }  \quad x \in [0, \infty)
\end{split}
\end{equation}

\begin{equation}
\begin{split}
G_{k}(x/p,\theta) &= \int_{0}^{U}\frac{p^{k}}{A(p,k)}U^{k-1}\frac{1}{1-p(1-U)^{k}}dU\\
                  &= \frac{p^{k}}{A(p,k)}\int_{0}^{U}\frac{U^{k-1}}{1-p(1-U)^{k}}dU\\
\end{split}
\end{equation}

let $y=1-p(1-U)$ ; $J=\frac{1}{p}$ ; $q=1-p$

\begin{equation}
\begin{split}
G_{k}(x/p,\theta) &= \frac{p^{k}}{A(p,k)}\int_{q}^{y}\big(\frac{y-q}{p}\big)^{k-1}\frac{1}{y^{k}}dy\\
                  &= \frac{1}{A(p,k)}\int_{q}^{y}(y-q)^{k-1}y^{-k}dy\\
                  &= \frac{1}{A(p,k)}\int_{q}^{y}\sum_{r=0}^{k-1} {k-1\choose r}(-q)^{r}y^{-(r+1)}dy\\
                  &= \frac{1}{A(p,k)}\sum_{r=0}^{k-1} {k-1\choose r}(-q)^{r}\int_{q}^{y}y^{-(r+1)}dy\\
                  &= \frac{1}{A(p,k)}\sum_{r=0}^{k-1} {k-1\choose r}(-q)^{r}\int_{q}^{y}I_{(x,r)}\\
                  &= \frac{1}{A(p,k)}\sum_{r=0}^{k-1} {k-1\choose r} (-1)^{r} (1-p)^{r}I_{(x,r)}
\end{split}
\end{equation}

where,
\begin{equation*}
    I_{(x,r)}=\int_{1-p}^{1-pe^{-\theta x}}t^{-(r+1)}dt
\end{equation*}

Otherwise,

Let $y=\frac{U}{1-p(1-U)}$ , then $U=\frac{y(1-p)}{1-py}$ and $1-p(1-U)=\frac{1-p}{1-py}$

$dU = \frac{1-p}{(1-py)^{2}}$

\begin{equation}
    G(x)=\int_{0}^{y}\frac{p^{k}y^{k-1}}{1-py} dy
\end{equation}

if $v=py$ ; $dv=pdy$

\begin{equation}
\begin{split}
    G(x)   &=\frac{1}{A(p,k)}\int_{0}^{v}\frac{v^{k-1}}{1-v} dv\\
           &=\frac{1}{A(p,k)}\sum_{r=0}^{\infty}\int_{0}^{v}v^{k+r-1} dv\\
           &=\frac{1}{A(p,k)}\sum_{r=0}^{\infty}\frac{v^{k+r}}{k+r}\\
           &=\frac{1}{A(p,k)}\sum_{z=k}^{\infty}\frac{v^{z}}{z} \mbox{ ; } z=k+r\\
           &=\frac{-ln(1-v)-\psi(k)\sum_{j=1}^{k-1}\frac{v^{j}}{j}}{A(p,k)}\\
           &=\frac{A(py,k)}{A(p,k)}
\end{split}
\end{equation}

\[
y=\frac{1-e^{-\theta x}}{1-pe^{-\theta x}}
\]

for $k=1$
\begin{equation}
    G(x)=\frac{ln\big(\frac{1-p}{1-pe^{-\theta x}}\big)}{ln(1-p)}
\end{equation}

The median is a solution of $G(x)=0.5$, then:

\begin{equation}
    x_\text{median}=-\frac{1}{\theta}\ln\big(\frac{1-\sqrt{1-p}}{p}\big)
\end{equation}

The mgf is $E(e^{t x})=E[(e^{-\theta x})^{-t/\theta}]$ ; let $u=e^{-\theta x}$

\begin{equation}\label{eq:generating-f-appendix}
\begin{split}
    E(e^{t x}) &= \frac{1}{A(p,k)}\int_{u=0}^{1}\frac{p^{k}u^{-t/\theta}(1-u)^{k-1}}{(1-pu)^{k}} du\\
               &= \frac{p^{k}}{A(p,k)}\int_{u=0}^{1}\sum_{i=0}^{\infty}{k-1+i\choose i}(pu)^{i}u^{-t/\theta}(1-u)^{k-1} du\\
               &= \frac{p^{k}}{A(p,k)}\sum_{i=0}^{\infty}{k-1+i\choose i}p^{i}\int_{u=0}^{1}u^{i-t/\theta}(1-u)^{k-1} du\\
               &=\frac{p^{k}}{A(p,k)}\sum_{i=0}^{\infty} {k-1+i\choose i} p^{i}\beta (i-\frac{t}{\theta}+1,k)
\end{split}
\end{equation}

where,
\begin{equation*}
    \beta (a,b)=\int_{0}^{1} t^{a-1} (1-t)^{b-1} dt
\end{equation*}

from the binomial theorem, we have:
\begin{equation*}
    (1-u)^{k-1}= \sum_{j=0}^{k-1}{k-1\choose j}(-u)^{j}=\sum_{j=0}^{k-1}{k-1\choose j}(-1)^{j} u^{j}
\end{equation*}

then,
\begin{equation}\label{eq:generating-f-appendix}
\begin{split}
    E(e^{t x}) &= \frac{p^{k}}{A(p,k)}\sum_{i=0}^{\infty}{k-1+i\choose i}p^{i}\int_{u=0}^{1}u^{i-t/\theta}(1-u)^{k-1} du\\
               &= \frac{p^{k}}{A(p,k)}\sum_{i=0}^{\infty}{k-1+i\choose i}p^{i}\sum_{j=0}^{k-1}{k-1\choose j}(-1)^{j}\int_{u=0}^{1}u^{i+j-t/\theta} du\\
               &= \frac{p^{k}}{A(p,k)}\sum_{i=0}^{\infty} \sum_{j=0}^{k-1}{k-1+i\choose i} {k-1\choose j}p^{i} (-1)^{j} \frac{1}{i+j-\frac{t}{\theta}+1}
\end{split}
\end{equation}

\begin{eqnarray*}
  E(x) &=& \frac{\partial}{\partial t} E(e^{t x})/_{t=0}\\
    &=& \frac{p^{k}}{\theta A(p,k)}\sum_{i=0}^{\infty} \sum_{j=0}^{k-1}{k-1+i\choose i} {k-1\choose j}p^{i} (-1)^{j} \frac{1}{(i+j-\frac{t}{\theta}+1)^2}/_{t=0} \\
    &=& \frac{p^{k}}{\theta A(p,k)}\sum_{i=0}^{\infty} \sum_{j=0}^{k-1}{k-1+i\choose i} {k-1\choose j}p^{i} (-1)^{j} \frac{1}{(i+j+1)^{2}}
\end{eqnarray*}

The $r^{th}$ moment is given by:

\begin{equation}\label{eq:rth-moment-appendix}
    E(x^{r})= \frac{\Gamma(r+1)}{\theta^{r}}\frac{p^{k}}{A(p,k)}\sum_{i=0}^{\infty} \sum_{j=0}^{k-1}{k-1+i\choose i} {k-1\choose j} p^{i} (-1)^{j} \frac{1}{(i+j+1)^{r+1}}
\end{equation}

Reliability or survival function:

\begin{equation}
\begin{split}
   S(x) &= Pr(X \geq x) = 1 - G(x) = \int_{x}^{\infty}f(t)d t\\
        &= 1-\frac{A(py,k)}{A(p,k)}\\
        &= 1-\frac{\ln(1-py)+\psi(k)\sum_{j=1}^{k-1}\frac{(py)^{j}}{j}}{\ln(1-p)+\psi(k)\sum_{j=1}^{k-1}\frac{p^{j}}{j}}\\
        &= 1-\frac{\ln\big(1-p\frac{1-e^{-\theta x}}{1-pe^{-\theta x}}\big)+\psi(k)\sum_{j=1}^{k-1}\frac{1}{j}
               \big(p\frac{1-e^{-\theta x}}{1-pe^{-\theta x}}\big)^{j}}{\ln(1-p)+\psi(k)\sum_{j=1}^{k-1}\frac{p^{j}}{j}}
\end{split}
\end{equation}

EM algorithm:

Using the joint distribution $g_{k}(z,x/p,\theta)$ we drive the conditional mass function as:

\begin{equation}\label{eq:prob-z-appendix}
   p(z/x,p,\theta)=\frac{\Gamma(z)}{\Gamma(k)\Gamma(z-k+1)}p^{z-k} e^{-\theta (z-k)x}(1-pe^{-\theta x})^{k}
\end{equation}

For $|a|<1$, we have

\begin{equation*}
    \sum_{i=0}^{\infty}{n+i\choose i}a^{i}=\frac{1}{(1-a)^{n+1}}
\end{equation*}

E-step:

\begin{equation}\label{eq:E-step-appendix}
\begin{split}
   E(z/x,p,\theta) &= \sum_{z=k}^{\infty} z p(z/x,p,\theta)\\
                   &= \sum_{z=k}^{\infty} k {z\choose k} (pe^{-\theta x})^{z-k}(1-pe^{-\theta x})^{k}\\
                   &= k (1-pe^{-\theta x})^{k}\sum_{z=k}^{\infty} {z\choose k} (pe^{-\theta x})^{z-k}\\
                   &= k (1-pe^{-\theta x})^{k}\sum_{z=k}^{\infty} {z\choose z-k} (pe^{-\theta x})^{z-k}\\
                   &= k (1-pe^{-\theta x})^{k}\sum_{t=0}^{\infty} {k+t\choose t} (pe^{-\theta x})^{t}\\
                   &= \frac{k}{1-pe^{-\theta x}}
\end{split}
\end{equation}

M-step:

Let $(x_{1}, x_{2},... ,x_{n})$ a random sample. The likelihood given the joint distribution $g_{k}(z,x/p,\theta)$ is:

\begin{equation}\label{eq:Lik-M-step-appendix}
\begin{split}
   L(p,\theta) &= \prod_{i=1}^{n}g_{k}(z,x/p,\theta)\\
               &= \prod_{i=1}^{n}\frac{\Gamma(z)}{\Gamma(k)\Gamma(z-k+1)}\frac{\theta p^{z} e^{-\theta (z-k+1) x} (1-e^{-\theta x})^{k-1}}{A(p,k)}\\
               &= \propto\prod_{i=1}^{n}\frac{\theta p^{z} e^{-\theta (z-k+1) x} (1-e^{-\theta x})^{k-1}}{A(p,k)}\\
               &= \propto \theta^{n} p^{\sum_{i=1}^{n}z_{i}} e^{-\theta\sum_{i=1}^{n}(z_{i}-k+1)x_{i}}A(p,k)^{-n}\prod_{i=1}^{n}(1-e^{-\theta x_{i}})^{k-1}\\
\end{split}
\end{equation}

The log-likelihood is:

\begin{equation}\label{eq:Lik-M-step-appendix}
       LL = n \ln(\theta)+\ln(p)\sum_{i=1}^{n}z_{i}-\theta\sum_{i=1}^{n}(z_{i}-k+1)x_{i}+\sum_{i=1}^{n}(1-e^{-\theta x_{i}})^{k-1}-n\ln(A(p,k))
\end{equation}

\begin{equation*}
    \frac{\partial LL}{\partial \theta} = \frac{n}{\theta}-\sum_{i=1}^{n}(z_{i}-k+1)x_{i}+(k-1)\sum_{i=1}^{n}\frac{x_{i}e^{-\theta x_{i}}}{1-e^{-\theta x_{i}}}=0
\end{equation*}

\begin{equation*}
\widehat{\theta}^{(r+1)}  =
n\bigg[\sum_{i=1}^{n}z_{i}x_{i}-(k-1)\sum_{i=1}^{n}\frac{x_{i}}{1-e^{-\theta^{(r+1)} x_{i}}}\bigg]^{-1}
\end{equation*}

we replace $z_{i}$ with  $E(z/x,p^{(r)},\theta^{(r)})=\frac{k}{1-p^{(r)}e^{-\theta^{(r)} x_{i}}}$

\begin{equation}\label{eq:M-step1-appendix}
\widehat{\theta}^{(r+1)}  =
n\bigg[\sum_{i=1}^{n}\frac{kx_{i}}{1-p^{(r)}e^{-\theta^{(r)} x_{i}}}-(k-1)\sum_{i=1}^{n}\frac{x_{i}}{1-e^{-\theta^{(r+1)} x_{i}}}\bigg]^{-1}
\end{equation}

\begin{equation*}
    \frac{\partial LL}{\partial p} = \frac{\sum_{i=1}^{n}z_{i}}{p}-n\frac{A(p,k)'}{A(p,k)}=0
\end{equation*}

where,
\begin{equation*}
   A(p,k) = \sum_{j=k}^{\infty}\frac{p^{j}}{j}
          = - \ln(1-p)-\psi(k)\sum_{j=1}^{k-1}\frac{p^{j}}{j}
\end{equation*}

\begin{equation*}
A(p,k)' =\frac{\partial A(p,k)}{\partial p} = \sum_{j=k}^{\infty}p^{j-1}=\frac{p^{k-1}}{1-p}
\end{equation*}

then,

\begin{equation*}
\frac{\sum_{i=1}^{n}z_{i}}{p}=n\frac{A(p,k)'}{A(p,k)}
\end{equation*}

\begin{equation}
\begin{split}
   p^{(r+1)} &= \frac{A(p,k)}{nA(p,k)'}\sum_{i=1}^{n}z_{i}\\
             &= \bigg[\frac{1-p^{(r+1)}}{n\Big(p^{(r+1)}\Big)^{k-1}}\sum_{j=k}^{\infty}\frac{\Big(p^{(r+1)}\Big)^{j}}{j}\bigg]\sum_{i=1}^{n}z_{i}\\
             &= \bigg[\frac{1-p^{(r+1)}}{n\Big(p^{(r+1)}\Big)^{k-1}}\sum_{j=k}^{\infty}\frac{\Big(p^{(r+1)}\Big)^{j}}{j}\bigg]\sum_{i=1}^{n}\frac{k}{1-p^{(r)}e^{-\theta^{(r)} x_{i}}}
\end{split}
\end{equation}

\begin{equation}\label{eq:M-step2-appendix}
\widehat{p}^{(r+1)}
  =
  \frac{-\big(1-p^{(r+1)}\big)\Big[\ln\big(1-p^{(r+1)}\big)+\psi(k)\sum_{j=1}^{k-1}\frac{\big(p^{(r+1)}\big)^{j}}{j}\Big]}{n\big(p^{(r+1)}\big)^{k-1}}\sum_{i=1}^{n}\frac{k}{1-p^{(r)}e^{-\theta^{(r)} x_{i}}}
\end{equation}

\newpage
\bibliographystyle{agsm}

\bibliography{lifetime}

\end{document}